\documentclass[graybox]{svmult}

\usepackage{type1cm}        
\usepackage{makeidx}         
\usepackage{graphicx}        
\usepackage{multicol}        
\usepackage[bottom]{footmisc}
\usepackage{comment}
\usepackage{multirow}
\usepackage{array}
\usepackage{mathtools}
\mathtoolsset{showonlyrefs}
\usepackage{booktabs}
\usepackage[percent]{overpic}
\usepackage{algorithm,algpseudocode}

\usepackage{newtxtext}       %
\usepackage[varvw]{newtxmath}       

\makeindex             

\newcommand{\rv}[1]{{\color{black} #1}}
\newcommand{\rvv}[1]{{\color{black} #1}}
\newcommand{\rvvv}[1]{{\color{black} #1}}

\newcommand{\numberset}{\mathbb}
\newcommand{\N}{\numberset{N}}
\newcommand{\R}{\numberset{R}}

\newcommand\edge{\varepsilon}

\newcommand{\LdueO}{L^2(\Omega)}

\newcommand{\HunoK}{H^1(\elementVEM)}
\newcommand{\HunozeroO}{H^1_0(\Omega)}
\newcommand{\elementVEM}{K}

\newcommand{\mesh}{\mathcal{T}}

\newcommand{\Puno}{\numberset{P}_1}

\newcommand{\PiNabla}{\Pi^\nabla}

\newcommand{\Acal}{\mathcal{A}}

\newcommand{\V}{V}

\newcommand{\Vrb}[1]{\V_{#1}^\mathrm{rb}}
\newcommand{\Wrb}[1]{W_{#1}^\mathrm{rb}}

\newcommand{\eiKrb}{ e_{M,i}^{\,\text{rb}}}

\newcommand{\vrb}{v^{\,\text{rb}}}

\newcommand{\wrb}{w^{\,\text{rb}}}

\renewcommand{\I}{\mathrm{Id}}

\DeclareFontFamily{U}{matha}{\hyphenchar\font45}
\DeclareFontShape{U}{matha}{m}{n}{
	<-6> matha5 <6-7> matha6 <7-8> matha7
	<8-9> matha8 <9-10> matha9
	<10-12> matha10 <12-> matha12
}{}
\DeclareSymbolFont{matha}{U}{matha}{m}{n}
\DeclareFontFamily{U}{mathx}{\hyphenchar\font45}
\DeclareFontShape{U}{mathx}{m}{n}{
	<-6> mathx5 <6-7> mathx6 <7-8> mathx7
	<8-9> mathx8 <9-10> mathx9
	<10-12> mathx10 <12-> mathx12
}{}
\DeclareSymbolFont{mathx}{U}{mathx}{m}{n}
\DeclareMathDelimiter{\vvvert} {0}{matha}{"7E}{mathx}{"17}%

\begin{document}

\title*{The Reduced Basis Multigrid scheme for the Virtual Element Method\thanks{This work was \rv{partially funded}  by MUR--PRIN/PNRR Bando 2022 (grant P2022BH5CB).
PFA is partially supported by the European Union (ERC SyG, NEMESIS, project number 101115663). Views and opinions expressed are, however, those of the authors only and do not necessarily reflect those of the European Union or the European Research Council Executive Agency. Neither the European Union nor the granting authority can be held responsible for them. The authors are members of INdAM--GNCS. }}
\author{Paola F. Antonietti, Silvia Bertoluzza, Fabio Credali}
\institute{Paola F. Antonietti, MOX, Dipartimento di Matematica, Politecnico di Milano, 20133 Milano, Italy, \email{paola.antonietti@polimi.it}\\
\and Silvia Bertoluzza, IMATI ``E. Magenes'', 27100 Pavia, Italy, \email{silvia.bertoluzza@imati.cnr.it}\\
\and Fabio Credali, King Abdullah University of Science and Technology, 23955 Thuwal, Saudi Arabia, \email{fabio.credali@kaust.edu.sa}}
%
%
\maketitle

\abstract{We present a non-nested W-cycle multigrid scheme for the lowest order Virtual Element Method on polygonal meshes. To avoid the implicit definition of the Virtual Element space, which poses several issues in the computation of intergrid operators that underpin multigrid methods, the proposed scheme uses a fully-conforming auxiliary space constructed by cheaply computing the virtual basis functions via the reduced basis method.}

\section*{Introduction}

The virtual element method (VEM) was introduced more than a decade ago as a generalization of the finite element method on general polytopal meshes~\cite{basic}. One of the main features of VEM is that the underlying discrete space is defined implicitly. Indeed, the basis functions are themselves solutions of a local differential problem in each mesh element. Such a local partial differential equation (PDE) is generally not solved explicitly, and the VEM discrete problem is defined using suitable polynomial projectors and stabilization terms. However, this construction presents some drawbacks: for instance, a wrong choice of the stabilization term may pollute the results~\cite{boffi}. From the algebraic solver perspective, the design of geometric multigrid schemes is non-trivial since the Virtual Element spaces are non-nested even if built on a sequence of nested meshes, and the injection operator cannot be defined by simple pointwise evaluation of basis functions at points located in the interior of the elements~\cite{multigrid}. To overcome the issues mentioned above, methods for reconstructing the virtual basis functions~\cite{credali,zerbinati,teora} have been introduced in the literature.

This work aims to exploit the reduced basis Virtual Element Method (\textsf{rb}VEM) we introduced in~\cite{credali} to construct a geometric W-cycle multigrid scheme for the lowest order virtual element method on non-nested meshes. To the best of our knowledge, this is the first work to exploit the efficient computation of the virtual basis function to design a multigrid algorithm for the $ h$-version of VEM.

We consider the standard virtual element formulation for the Poisson problem and the \textsf{rb}VEM space as an auxiliary tool for constructing the intergrid operators between mesh levels. In particular, we define the prolongation operator as the $L^2-$projection~\cite{dg} between \textsf{rb}VEM spaces on different levels, without the need of modifying the original discrete problem.  We present numerical tests showing that the convergence factor of our new scheme is independent of the resolution of the underlying fine mesh. We also show that our algorithm is independent of the number of reduced basis functions employed for constructing the auxiliary \textsf{rb}VEM space. A forthcoming work will address a detailed convergence analysis for the proposed method, as well as the design of V-cycle schemes.

\section{Model problem and Virtual Element discretization}

We consider the weak formulation of the Poisson equation with homogeneous Dirichlet boundary conditions in a two-dimensional polygonal domain $\Omega\subset\R^2$, which reads: $\text{given $f\in\LdueO$, find $u\in \HunozeroO$ such that}$
\begin{equation}\label{eq:poisson}
    \begin{aligned}
 &\Acal(u,v) = (\nabla u,\nabla v)_\Omega = (f,v)_\Omega\qquad\forall v\in\HunozeroO,
    \end{aligned}
\end{equation}
where the notation $(\cdot,\cdot)_\Omega$ indicates the scalar product in $\LdueO$.

We discretize our problem by considering the virtual element method. Let $\mesh_h$ be a generic polygonal tessellation of $\Omega$ made of disjoint star-shaped polygonal elements $\elementVEM$ with diameter $h_\elementVEM$. The size of $\mesh_h$ is denoted by $h=\max_{\elementVEM\in\mesh_h}h_\elementVEM$. The local virtual element space~\cite{basic} of degree one  is defined as $\V(\elementVEM) = \{v\in\HunoK:\,v_{|\edge}\in\Puno(\edge)$ $\forall\edge\in\partial\elementVEM,\,\Delta v=0\,\text{in }\elementVEM\}$,
and the global discrete space is then obtained by gluing all local spaces by continuity
\begin{equation}
    \V(\mesh_h) = \{v\in\HunozeroO:\,v\in\V(\elementVEM)\,\forall\elementVEM\in\mesh_h\}.
\end{equation}
We observe that a function $v\in\V(\elementVEM)$ is uniquely identified by its values at the $N$ vertices of $\elementVEM$, which are a unisolvent set of degrees of freedom. Moreover, accuracy is guaranteed since $\Puno(\elementVEM)\subset\V(\elementVEM)$.

In order to solve Equation~\eqref{eq:poisson} by means of a geometric multigrid algorithm, we consider a sequence of quasi-uniform tessellations $\{\mesh_j\}_{j=1}^J$ with mesh sizes $h_j$ satisfying the constraint 
    $C\,h_{j-1} \le h_j \le h_{j-1}$,
for a positive constant $C$. Moreover, $\{\mesh_j\}_{j=1}^J$ are not-nested. We do not assume that meshes are nested, as the nestedness of the meshes does not imply that of the corresponding spaces. We then seek for the virtual element approximation of $u$ at the finest level grid $J$:
\begin{equation}\label{eq:poisson_h}
    \begin{aligned}
        &\text{find $u_J\in \V(\mesh_J)$ such that}
        \quad
        \Acal_J(\rv{u_J},v_J) = (f_J,v_J)_\Omega\qquad\forall v_J\in\V(\mesh_J).
    \end{aligned}
\end{equation}
The discrete bilinear form $\Acal_J$ is defined as 
\begin{equation}
    \begin{aligned}
    &\Acal_J(u_J,v_J) = \sum_{\elementVEM_J\in\mesh_J} \Acal_J^{\elementVEM_J}(u_J,v_J)\rv{,}\\
    &\Acal_J^{\elementVEM_J}(u_J,v_J)=(\nabla\PiNabla u_J,\nabla\PiNabla v_J)_{\elementVEM_J} + S_J^{\elementVEM_J}((\I-\PiNabla)u_J,(\I-\PiNabla)v_J),
    \end{aligned}
\end{equation}
where $\PiNabla$ is the elliptic projection operator onto $\Puno(\elementVEM_J)$ solving
\begin{equation}
    (\nabla v_J-\nabla\PiNabla v_J,\nabla q)_{\elementVEM_J} = 0 \quad \forall q\in\Puno(\elementVEM_J),
    \quad
    \int_{\partial\elementVEM_J} \PiNabla v_J\,ds = \int_{\partial\elementVEM_J} v_J\,ds,
\end{equation}
and $S_J^{\elementVEM_J}$ is the semi-positive definite bilinear form defined as
\begin{equation}\label{eq:stab}
    S_J^{\elementVEM_J}(w_J,v_J) = \sum_{i=1}^N w_J(\mathbf{x}_i)v_J(\mathbf{x}_i),
\end{equation}
with $\mathbf{x}_1,\dots,\mathbf{x}_N$ being the vertices of $\elementVEM_J\in\mesh_J$. 
\rv{Finally, the discrete right-hand side is defined as
$(f_J,v_J)_\Omega = \sum_{\elementVEM_J\in\mesh_J} (\Pi^0_0 f,v_J)_{\elementVEM_J}$, where $\Pi^0_0$ is the $L^2$ projection onto constants, i.e.,  $\Pi^0_0v_J=|\elementVEM_J|^{-1}\int_{\elementVEM_J}v_J\,dx$ for any $v_J\in\V(\elementVEM_J)$.
}
Eq.~\eqref{eq:poisson_h} can be reformulated in terms of operators as $ A_Ju_J=f_J$.
\section{The reduced basis multigrid algorithm}
In this section we design a W-cycle non-nested multigrid algorithm to solve Eq.~\eqref{eq:poisson_h}. A key ingredient for the construction of the method are the intergrid operators. Since $\V(\mesh_{j-1})\not\subset\V(\mesh_j)$, a good choice of prolongation operator is the $L^2$ projection~\cite{dg} $I_{j-1}^j:\V(\mesh_{j-1})\rightarrow\V(\mesh_j)$ defined, for $v_{j-1}\in\V(\mesh_{j-1})$, as
\begin{equation}
    (I_{j-1}^j v_{j-1},w_j)_\Omega = (v_{j-1},w_j)_\Omega\qquad\forall w_j\in\V(\mesh_j).
\end{equation}
Clearly, due to the implicit definition of the virtual element spaces, the quantities above are not fully computable through the degrees of freedom. We thus consider the \textsf{rb}VEM space introduced in~\cite{credali}, whose basis functions are explicitly constructed and have the same degrees of freedom as the original VEM space. For a generic element $\elementVEM_j\in\mesh_j$, the \textsf{rb}VEM is defined as
\begin{equation}\label{eq:rbvem_space}
    \Vrb{}(\elementVEM_j) = \Puno(\elementVEM_j)\oplus\Wrb{}(\elementVEM_j),
\end{equation}
where $\Wrb{}(\elementVEM_j)\subset{\textrm{span}}\{\eiKrb,\,i=1,\dots,N\}$ is such that $\wrb_j=\sum_{i=1}^N \mathsf{w}_i\eiKrb\in\Wrb{}(\elementVEM_j)$ if and only if $\PiNabla(\sum_{i=1}^N \mathsf{w}_i e_i)=0$, with $\{e_i\}_{i=1}^N$ being the shape functions of $\V(\elementVEM_j)$.
The functions $\eiKrb$ are an approximation of the virtual element basis functions computed by means of the reduced basis method with complexity, i.e. the dimension of the reduced basis space, equal to $M$. \rvv{By construction, $\eiKrb(\mathbf{x}_k)=\delta_{i,k}$ for $i,k=1,\dots,N$}, \rvvv{see~\cite{credali}}.  In addition, the space $\Wrb{}(\elementVEM_j)$ is constructed to be $\Acal-$orthogonal to $\Puno(\elementVEM_j)$. Given $v_j\in\V(\elementVEM_j)$, its counterpart $\vrb_j\in\Vrb{}(\elementVEM_j)$ is constructed as $\vrb_j=\PiNabla v_j+v_j^\perp$, where $v_j^\perp\in\Wrb{}(\elementVEM_j)$ is the linear combination $v_j^\perp=\sum_{i=1}^N (v_j(\mathbf{x}_i)-\PiNabla v_j(\mathbf{x}_i))\eiKrb$. Given $v_j,w_j\in\V(\elementVEM_j)$ and their counterparts $\wrb_j,\vrb_j\in\Vrb{}(\elementVEM_j)$, we set
$$
\begin{aligned}
\Acal_j^{\elementVEM_j}&(\wrb_j,\vrb_j)
= (\nabla\PiNabla w_j,\nabla\PiNabla v_j)_{\elementVEM_j} + S_j^{\elementVEM_j}(w_j^\perp,v_j^\perp)\\
&\rvv{=(\nabla\PiNabla w_j,\nabla\PiNabla v_j)_{\elementVEM_j} + \sum_{i=1}^N w_j^\perp(\mathbf{x}_i)v_j^\perp(\mathbf{x}_i)}\\
&\rv{= (\nabla\PiNabla w_j,\nabla\PiNabla v_j)_{\elementVEM_j} + \sum_{i=1}^N
\big[(w_j-\PiNabla w_j)(\mathbf{x}_i)\big]\big[(v_j-\PiNabla v_j)(\mathbf{x}_i)\big]}\\
&= (\nabla\PiNabla w_j,\nabla\PiNabla v_j)_{\elementVEM_j} + S_j^{\elementVEM_j}((\I-\PiNabla)w_j,(\I-\PiNabla)v_j)
=\Acal_j^{\elementVEM_j}(w_j,v_j),
\end{aligned}
$$
\rvv{where we exploited that, for $k=1,\dots,N$, 
$$
\aligned
v_j^\perp(\mathbf{x}_k)
&=\sum_{i=1}^N [v_j(\mathbf{x}_i)-\PiNabla v_j(\mathbf{x}_i)]\eiKrb(\mathbf{x}_k)
=\sum_{i=1}^N [v_j(\mathbf{x}_i)-\PiNabla v_j(\mathbf{x}_i)]\delta_{i,k}\\
&= v_j(\mathbf{x}_k)-\PiNabla v_j(\mathbf{x}_k)
\endaligned
$$
and, analogously, $w_j^\perp(\mathbf{x}_k)=w_j(\mathbf{x}_k)-\PiNabla w_j(\mathbf{x}_k)$.
}
\rv{As for the standard VEM space, the global \textsf{rb}VEM space is obtained by gluing all local spaces by continuity, that is
\begin{equation}
    \Vrb{}(\mesh_j) = \{\vrb_j\in\HunozeroO:\,\vrb_j\in\Vrb{}(\elementVEM_j)\,\forall\elementVEM_j\in\mesh_j\},\qquad j=1,\dots,J.
\end{equation}}
We then solve Eq.~\eqref{eq:poisson_h} in the original VEM space, by employing the newly introduced \textsf{rb}VEM as a support for computing intergrid operators. Given $\vrb_{j-1}\in\Vrb{}(\mesh_{j-1})$, we define the prolongation operator $I_{j-1}^{j}:\Vrb{}(\mesh_{j-1})\rightarrow\Vrb{}(\mesh_{j})$ as
\begin{equation}\label{eq:prolongation}
    (I_{j-1}^j \vrb_{j-1},\wrb_j)_\Omega = (\vrb_{j-1},\wrb_j)_\Omega\qquad\forall \wrb_j\in\Vrb{}(\mesh_j),
\end{equation}
while the restriction operator $I_{j}^{j-1}:\Vrb{}(\mesh_{j})\rightarrow\Vrb{}(\mesh_{j-1})$ is defined as its adjoint
\begin{equation}\label{eq:restriction}
    (I_{j}^{j-1} \wrb_j,\vrb_{j-1})_\Omega = (\wrb_j,\vrb_{j-1})_\Omega\qquad\forall \vrb_{j-1}\in\Vrb{}(\mesh_{j-1}).
\end{equation}
In our scheme, the prolongation and restriction operators between the VEM spaces $\V(\mesh_{j-1})$, $\V(\mesh_j)$ inherit their definition from~\eqref{eq:prolongation} and~\eqref{eq:restriction}, respectively.
\begin{algorithm}[t!]
    \caption{Multigrid W-cycle iteration}
    \label{alg:iteration}
    \begin{algorithmic}
    \State Initialize $u_0\in\V(\mesh_J)$:
    \For{$k=0,1,\dots$}
    \State $u_{k+1}=\mathsf{MG}_\mathcal{W}(J,f_J,u_k,m_1,m_2);$
    \State $u_k = u_{k+1};$
    \EndFor
    \end{algorithmic}
\end{algorithm}
We point out that the scalar product on the \rv{left-hand side of~\eqref{eq:prolongation}}, which involves functions in $\Vrb{}(\mesh_j)$, can be efficiently computed by exploiting the properties of the reduced basis method~\cite{credali}. On the other hand, the scalar product at the \rv{right-hand side involves} functions defined on two different mesh levels and requires the computation of the intersection between $\mesh_j$ and $\mesh_{j-1}$, as commented in~\cite{dg}. Two meshes can be efficiently intersected if a bounding box technique is employed to detect cell-by-cell interactions.
As a smoothing scheme, we choose the Richardson iteration. Given $\V^\star(\mesh_j)$ the dual space of $\V(\mesh_j)$, let $A_j:\V(\mesh_j)\rightarrow\V^\star(\mesh_j)$ be the operator defined as
$(A_jw,v)_\Omega = \Acal_j(w,v)$ for all $w,v\in\V(\mesh_j)$
and let $\Lambda_j\in\R$ denote the maximum eigenvalue of $A_j$, for $j=2,\dots,J$. The smoothing scheme is then defined as $B_j = \Lambda_j \I_j$, $j=2,\dots,J$.
\begin{algorithm}[htbp!]
\caption{Multigrid W-cycle scheme: one iteration at level \textcolor{black}{$j\ge1$}}
\label{alg:multilevel}
\begin{algorithmic}
\If{$j=1$}
\State \textcolor{black}{$\mathsf{MG}_\mathcal{W} (1,g,z_0,m_1,m_2)=A_1^{-1}g.$}
\Else
\For{$i=1,\dots,m_1$}
\hfill{{\it (Pre-smoothing)}}
\State $z^{(i)}=z^{(i-1)}+B_j^{-1}(g-A_jz^{(i-1)});$ \EndFor\vspace{0.2cm}
\State $r_{j-1} = I_j^{j-1}(g-A_jz^{(m_1)})$;
\hfill{{\it (Coarse grid correction)}}
\State $\overline{e}_{j-1} = \mathsf{MG}_\mathcal{W} (j-1,r_{j-1},0,m_1,m_2)$;
\State $e_{j-1} = \mathsf{MG}_\mathcal{W} (j-1,r_{j-1},\overline{e}_{j-1},m_1,m_2)$;
\State $z^{(m_1+1)}=z^{(m_1)}+I_{j-1}^je_{j-1}$;\vspace{0.2cm}
\For{$i=m_1+2,\dots,m_1+m_2+1$} 
\hfill{{\it (Post-smoothing)}}\State $z^{(i)}=z^{(i-1)}+B_j^{-1}(g-A_jz^{(i-1)});$ \EndFor\vspace{0.3cm}
\State $\mathsf{MG}_\mathcal{W} (j,g,z_0,m_1,m_2)=z^{(m_1+m_2+1)}.$
\EndIf
\end{algorithmic}
\end{algorithm}
The W-cycle iteration for computing $u_J$ is summarized in Algorithm~\ref{alg:iteration}, where $u_0\in\V(\mesh_J)$ denotes a suitable initial guess, and where $m_1,m_2\in\N$ denote the number of pre- and post-smoothing steps, respectively. We observe that $\mathsf{MG}_\mathcal{W}(J,f_J,u_k,m_1,m_2)$ represents the approximate solution obtained after one iteration of the non-nested W-cycle, which is defined by the recursive procedure described in Algorithm~\ref{alg:multilevel}. Indeed, the W-cycle is defined by induction: given the general problem of finding $z\in\V(\mesh_j)$ satisfying
\begin{equation}\label{eq:pbz}
A_j z=g
\end{equation}
for $j\ge2$ and $g\in\LdueO$, then $\mathsf{MG}_\mathcal{W}(j,g,z_0,m_1,m_2)$ is the approximate solution of~\eqref{eq:pbz} obtained after one multigrid iteration with initial guess $z_0\in\V(\mesh_j)$ and $m_1,m_2$ smoothing steps. Clearly, on the coarsest level $j=1$, the problem is solved exactly.

\begin{remark}
    In this work, we are using the spaces $\Vrb{}(\mesh_{j})$, $j=1,\dots,J$ as an auxiliary tool for computing intergrid operators. At the same time, the discrete problem is formulated in the standard virtual element space $\V(\mesh_J)$. We observe that, if the reduced basis stabilization introduced in~\cite{credali} is considered, then Eq.~\eqref{eq:poisson_h} can be reformulated in $\Vrb{}(\mesh_{J})$, so that $\Acal_J(w_J^\mathrm{rb},v_J^\mathrm{rb})=\Acal(w_J^\mathrm{rb},v_J^\mathrm{rb})$ for $w_J^\mathrm{rb},v_J^\mathrm{rb}\in\Vrb{}(\mesh_J)$. We thus obtain a discrete formulation in the fully-conforming \textsf{rb}VEM space, and our multigrid scheme remains unchanged.
    \rv{As discussed in~\cite{credali}, the use of \textsf{rb}VEM as a discretization method would require more computational effort than computing the stabilization term in~\eqref{eq:stab}, which is extremely cheap. However, this overhead is not worthwhile for isotropic problems, as VEM and \textsf{rb}VEM exhibit the same approximation properties. Conversely, \textsf{rb}VEM is more accurate than standard VEM when solving anisotropic problems, which are out of the scope of the present work.   
    We finally note that the proposed Reduced Basis Multigrid scheme is non-intrusive, as it does not require any modifications to existing VEM and Multigrid implementations.
    }
\end{remark}

\section{Numerical results}
\begin{figure}[!htbp]
    \centering
    \includegraphics[width=0.7\linewidth]{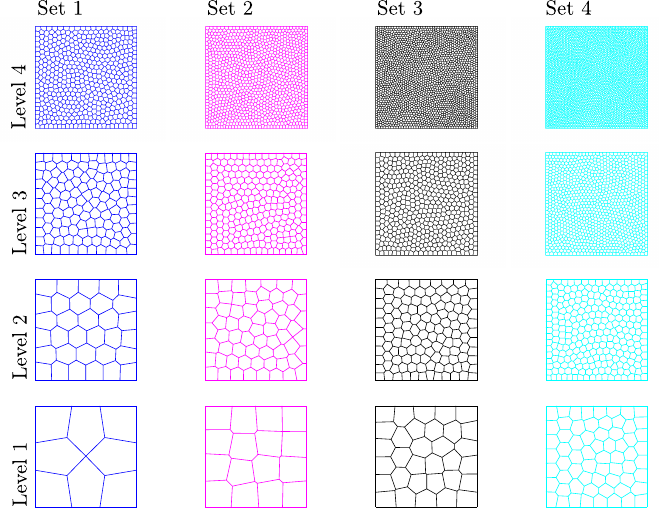}
    \caption{Sets of non-nested grids employed for numerical tests.}
    \label{fig:meshes}
\end{figure}

In this section, we present numerical tests to assess the performance of the proposed multigrid scheme. We set $\Omega=(0,1)^2$, and we consider the sets of non-nested meshes shown in Figure~\ref{fig:meshes}. The finest meshes (Level~4) are made of 512 (Set~1), 1024 (Set~2), 2048 (Set~3), and 4096 (Set~4) elements. We remark that each mesh was generated independently of the others. Moreover, the ratio between the numbers of elements of a coarse mesh and of the corresponding fine mesh is equal to~1/4. In the following, we report the iteration counts needed to reduce the relative residual below the tolerance $tol=10^{-8}$. We also compute the convergence factor
\begin{equation}
    \rho_J = \exp\left(\frac1{n_{it,J}}\log\frac{\|\mathbf{r}_{n_{it,J}}\|}{\|\mathbf{r}_0\|}\right),
\end{equation}
where $n_{it,J}$ is the iteration count for reducing the error below the considered tolerance by means of the W-cycle scheme with $J$ levels ($J=2,3,4$), while $\mathbf{r}_{n_{it,J}}$ and $\mathbf{r}_0$ are the final and initial residual vectors, respectively. In all our tests, we set $m_1=m_2=m$.
We first assess the independence of the multigrid scheme from the number $M$ of reduced basis functions employed to construct the shape functions of the space $\Wrb{}(\elementVEM_j)$ in~\eqref{eq:rbvem_space}. By looking at Table~\ref{tab:w_rich_m}, it is clear that the value of $M$ does not affect the convergence of the proposed scheme: indeed, for fixed $m$ and for a fixed number of levels $J$, the iteration count is the same for $M=1$ and $M=50$. For the next test, we thus set $M=1$.
\begin{table}[t!] 
\caption{Iterations count of the W-Cycle multigrid as a function of $m$. Comparison in terms of the number of reduced basis functions used to build the auxiliary space: $M=1$ vs $M=50$.}
\setlength\extrarowheight{1pt}
\centering
\begin{tabular}{p{1.25cm}>{\centering\arraybackslash}p{0.725cm}>{\centering\arraybackslash}p{0.725cm}>{\centering\arraybackslash}p{0.725cm}>{\centering\arraybackslash}p{0.725cm}>{\centering\arraybackslash}p{0.725cm}>{\centering\arraybackslash}p{0.725cm}p{0.25cm}>{\centering\arraybackslash}p{0.725cm}>{\centering\arraybackslash}p{0.725cm}>{\centering\arraybackslash}p{0.725cm}>{\centering\arraybackslash}p{0.725cm}>{\centering\arraybackslash}p{0.725cm}>{\centering\arraybackslash}p{0.725cm}}
\toprule[1.25pt]
\multicolumn{14}{c}{\textbf{W-cycle}, Richardson smoother, $tol=10^{-8}$}\\
\toprule[1.25pt]
 & \multicolumn{6}{c}{Set 1} & & \multicolumn{6}{c}{Set 2}\\ 
 \cline{2-7}\cline{9-14}
& \multicolumn{2}{c}{2 levels} & \multicolumn{2}{c}{3 levels} & \multicolumn{2}{c}{4 levels} & &\multicolumn{2}{c}{2 levels} &  \multicolumn{2}{c}{3 levels} & \multicolumn{2}{c}{4 levels}\\
M& $1$ & $50$ & $1$ & $50$ & $1$ & $50$ &  & $1$ & $50$ & $1$ & $50$ & $1$ & $50$ \\
 \hline
$m=3$ & 8&8 & 9&9 & 8&8 & & 8&8 & 8&8 & 8&8\\
$m=6$ & 7&7 & 7&7 & 6&6 & & 6&6 & 6&6 & 6&6\\
$m=8$ & 6&6 & 5&5 & 5&5 & & 6&6 & 5&5 & 5&5\\
\toprule[1.25pt]
 & \multicolumn{6}{c}{Set 3} & & \multicolumn{6}{c}{Set 4}\\ 
 \cline{2-7}\cline{9-14}
& \multicolumn{2}{c}{2 levels} & \multicolumn{2}{c}{3 levels} & \multicolumn{2}{c}{4 levels} & &\multicolumn{2}{c}{2 levels} &  \multicolumn{2}{c}{3 levels} & \multicolumn{2}{c}{4 levels}\\
M& $1$ & $50$ & $1$ & $50$ & $1$ & $50$ &  & $1$ & $50$ & $1$ & $50$ & $1$ & $50$ \\
 \hline
$m=3$ & 8&8 & 8&8 & 8&8 & & 8&8 & 8&8 & 8&8\\
$m=6$ & 6&6 & 6&6 & 6&6 & & 6&6 & 6&6 & 6&6\\
$m=8$ & 5&5 & 5&5 & 5&5 & & 5&5 & 5&5 & 5&5\\
\toprule[1.25pt]
\end{tabular}
\label{tab:w_rich_m}
\end{table}
The second test studies the convergence properties of our W-cycle scheme. In Table~\ref{tab:w_rich_conv}, we report both the convergence factor and the iteration count for all sets of meshes by varying the number $m$ of pre- and post-smoothing steps. It is clear that the convergence factor does not depend on the mesh size $h$. Indeed, for fixed $J=2,3,4$ and $m=3,6,8$, the convergence factor is roughly constant. This fact implies that the iteration counts needed by the W-cycle to reduce the residual below the given tolerance are independent of the resolution of the underlying mesh. Moreover, as $m$ increases, the convergence factor decreases.
\begin{table}[t!]
\caption{Convergence factor $\rho$ and iterations count of the W-Cycle multigrid as a function of $m$. The auxiliary space is constructed using a single reduced basis function.}
\setlength\extrarowheight{2pt}
\centering
\resizebox{\textwidth}{!}{
\begin{tabular}{p{1.25cm}p{1.45cm}p{1.45cm}p{1.45cm}p{0.25cm}p{1.45cm}p{1.45cm}p{1.45cm}p{1.45cm}}
\toprule[1.25pt]
\multicolumn{8}{c}{\textbf{W-cycle}, Richardson smoother, $tol=10^{-8}$, $M=1$}\\
\toprule[1.25pt]
 & \multicolumn{3}{c}{Set 1} & & \multicolumn{3}{c}{Set 2}\\ 
 \cline{2-4}\cline{6-8}
 & 2 levels & 3 levels & 4 levels & &2 levels &  3 levels & 4 levels\\
 \hline
$m=3$ & 0.0958 (8) & 0.1174 (9) & 0.0964 (8) & & 0.0843 (8) & 0.0896 (8) & 0.0991 (8)\\
$m=6$ & 0.0564 (7) & 0.0515 (7) & 0.0408 (6) & & 0.0428 (6) & 0.0364 (6) & 0.0392 (6)\\
$m=8$ & 0.0306 (6) & 0.0246 (5) & 0.0217 (5) & & 0.0310 (6) & 0.0204 (5) & 0.0199 (5)\\
\toprule[1.25pt]
 & \multicolumn{3}{c}{Set 3} & & \multicolumn{3}{c}{Set 4}\\ 
\cline{2-4}\cline{6-8}
 & 2 levels & 3 levels & 4 levels & & 2 levels &  3 levels & 4 levels\\
\hline
$m=3$ & 0.0805 (8) & 0.0825 (8) & 0.0942 (8) & & 0.0977 (8) & 0.0796 (8) & 0.0984 (8)\\
$m=6$ & 0.0369 (6) & 0.0308 (6) & 0.0384 (6) & & 0.0418 (6) & 0.0323 (6) & 0.0397 (6)\\
$m=8$ & 0.0229 (5) & 0.0174 (5) & 0.0196 (5) & & 0.0211 (5) & 0.0128 (5) & 0.0188 (5)\\
\toprule[1.25pt]
\end{tabular}
}
\label{tab:w_rich_conv}
\end{table}

\section{Conclusions}

We introduced a non-nested W-cycle multigrid scheme for the lowest order virtual element method based on the \textsf{rb}VEM auxiliary space introduced in~\cite{credali}. Numerical tests showed that the W-cycle multigrid scheme is scalable and that its convergence factor does not depend on the number of reduced basis functions used to construct the auxiliary space. A forthcoming paper will present an exhaustive convergence analysis of the proposed method, as well as extensions to the V-cycle algorithm and the application of smoothing schemes other than the Richardson iteration.

\bibliographystyle{abbrv}
\bibliography{biblio}

@article{basic,
	title={Basic principles of virtual element methods},
	author={Beir{\~a}o da Veiga, Lourenco and Brezzi, Franco and Cangiani, Andrea and Manzini, Gianmarco and Marini, L Donatella and Russo, Alessandro},
	journal={Mathematical Models and Methods in Applied Sciences},
	volume={23},
	number={01},
	pages={199--214},
	year={2013},
	publisher={World Scientific}
}

@article{boffi,
	title={Approximation of {PDE} eigenvalue problems involving parameter dependent matrices},
	author={Boffi, Daniele and Gardini, Francesca and Gastaldi, Lucia},
	journal={Calcolo},
	volume={57},
	number={4},
	pages={41},
	year={2020},
	publisher={Springer}
}

@article{multigrid,
	title={Agglomeration-based geometric multigrid schemes for the virtual element method},
	author={Antonietti, Paola F and Berrone, Stefano and Busetto, Martina and Verani, Marco},
	journal={SIAM Journal on Numerical Analysis},
	volume={61},
	number={1},
	pages={223--249},
	year={2023},
	publisher={SIAM}
}

@article{credali,
	title={Reduced basis stabilization and post-processing for the virtual element method},
	author={Credali, Fabio and Bertoluzza, Silvia and Prada, Daniele},
	journal={Computer Methods in Applied Mechanics and Engineering},
	volume={420},
	pages={116693},
	year={2024},
	publisher={Elsevier}
}

@article{zerbinati,
	title={When rational functions meet virtual elements: the lightning virtual element method},
	author={Trezzi, Manuel and Zerbinati, Umberto},
	journal={Calcolo},
	volume={61},
	number={3},
	pages={35},
	year={2024},
	publisher={Springer}
}

@article{teora,
	title={The lowest-order neural approximated virtual element method on polygonal elements},
	author={Berrone, Stefano and Pintore, Moreno and Teora, Gioana},
	journal={Computers \& Structures},
	volume={314},
	pages={107753},
	year={2025},
	publisher={Elsevier}
}

@article{dg,
	title={V-cycle multigrid algorithms for discontinuous {G}alerkin methods on non-nested polytopic meshes},
	author={Antonietti, Paola F and Pennesi, Giorgio},
	journal={Journal of Scientific Computing},
	volume={78},
	number={1},
	pages={625--652},
	year={2019},
	publisher={Springer}
}
\end{document}